\documentclass{amsart} 

\oddsidemargin \evensidemargin
\newtheorem{theorem}{Theorem}[section]
\newtheorem{proposition}{Proposition}[section]

\newtheorem{lemma}[theorem]{Lemma}
\newtheorem{cor}[theorem]{Corollary}
\theoremstyle{definition}

\theoremstyle{remark}

\numberwithin{equation}{section}
\numberwithin{figure}{section}
\usepackage{color}
\usepackage{amsmath, latexsym, amssymb}
\usepackage{graphicx}
\usepackage{amsfonts}
\usepackage{subfigure}
\usepackage{amsthm}
\usepackage{amsmath}
\usepackage{graphics}
\usepackage{amssymb}
\usepackage{latexsym}
\usepackage{amscd}
\usepackage{color}
\usepackage{graphicx}
\usepackage{subfigure}
\usepackage{pict2e}
\usepackage[all]{xy}
\usepackage{setspace}

\DeclareMathOperator{\col}{col}
\DeclareMathOperator{\key}{key}
\DeclareMathOperator{\colform}{colform}
\DeclareMathOperator{\st}{\mathfrak{U}(\omega,\lambda)}

\newlength\cellsize 
\savebox2{%
\begin{picture}(15,15)
\put(0,0){\line(1,0){15}}
\put(0,0){\line(0,1){15}}
\put(15,0){\line(0,1){15}}
\put(0,15){\line(1,0){15}}
\end{picture}}
\newcommand\cellify[1]{\def\thearg{#1}\def\nothing{}%
\ifx\thearg\nothing
\vrule width0pt height\cellsize depth0pt\else
\hbox to 0pt{\usebox2\hss}\fi%
\vbox to 15\unitlength{
\vss
\hbox to 15\unitlength{\hss$#1$\hss}
\vss}}
\newcommand\tableau[1]{\vtop{\let\\=\cr
\setlength\baselineskip{-16000pt}
\setlength\lineskiplimit{16000pt}
\setlength\lineskip{0pt}
\halign{&\cellify{##}\cr#1\crcr}}}
\savebox3{%
\begin{picture}(15,15)
\put(0,0){\line(1,0){15}}
\put(0,0){\line(0,1){15}}
\put(15,0){\line(0,1){15}}
\put(0,15){\line(1,0){15}}
\end{picture}}
\newcommand\expath[1]{%
\hbox to 0pt{\usebox3\hss}%
\vbox to 15\unitlength{
\vss
\hbox to 15\unitlength{\hss$#1$\hss}
\vss}}

\author[S. Mason]{S. Mason}
\thanks{Partially supported by NSF postdoctoral research fellowship 
               DMS-0603351 (S.M.)}
\address{Department of Mathematics, Davidson College}

\email{samason@davidson.edu}
\urladdr{http://www.davidson.edu/math/mason}

\title[Demazure atoms]{An explicit construction of type A Demazure atoms}



\begin{document}

\begin{abstract}

Demazure characters of type A, which are equivalent to key polynomials, have been decomposed by Lascoux and Sch\"{u}tzenberger into standard bases.  We prove that the resulting polynomials, which we call Demazure atoms, can be obtained from a certain specialization of nonsymmetric Macdonald polynomials.  This combinatorial interpretation for Demazure atoms accelerates the computation of the right key associated to a semi-standard Young tableau.  Utilizing a related construction, we provide a new combinatorial description of the key polynomials.\\

\end{abstract}

\maketitle

\section{Introduction}

The Demazure character formula generalizes the Weyl character formula to highest-weight modules over symmetrizable Kac-Moody Lie algebras.  In particular, if $V(\lambda)$ is a highest weight module of weight $\lambda$, then the extremal weight vector $u_{\omega \lambda}$ of weight $\omega \lambda$ generates a $U(\mathfrak{n})$-submodule $U(\mathfrak{n})u_{\omega \lambda}$.  The formal character of this submodule is given by Demazure's character formula~\cite{D1},~\cite{Jo85}.  The Demazure characters corresponding to the general linear Lie algebra $\mathfrak{g}l_n(\mathbb{C})$ are equivalent to the key polynomials, which are described~\cite{RS95} as the sums of the weights of semi-standard Young tableaux (SSYT) whose right key is bounded by a certain key $K(\omega,\lambda)$.  

Lascoux and Sch\"{u}tzenberger~\cite{LS90} study the smallest non-intersecting pieces, $\mathfrak{U}(\omega,\lambda)$, of type $A$ Demazure characters.  They call the resulting polynomials {\it standard bases} and describe them combinatorially as the sums of the weights of all semi-standard Young tableaux whose right key is equal to the key $K(\omega, \lambda)$.  Each semi-standard Young tableau appears in precisely one such polynomial, implying that the polynomials $\st$ form a decomposition of the Schur functions.

There exists a decomposition of the Schur functions into the polynomials $E_{\gamma}(x;0,0)$, which are obtained by setting $q=t=0$ in the combinatorial formula for integral form nonsymmetric Macdonald polynomials~\cite{HHL05c}.  The $E_{\gamma}(x;0,0)$ are obtained from the weights of {\it semi-skyline augmented fillings}, which are fillings of composition diagrams with positive integers in such a way that the columns are weakly decreasing and the rows satisfy an inversion condition.  Semi-skyline augmented fillings are in bijection with semi-standard Young tableaux and satisfy a variation of the Robinson-Schensted-Knuth algorithm~\cite{M06}.

\begin{theorem}{\label{atoms}}
The standard base $\mathfrak{U}(\omega, \lambda)$ is equal to the specialized nonsymmetric Macdonald polynomial $E_{\omega(\lambda)} (X; 0,0)$.
\end{theorem}

We obtain an efficient method for computing the right key of a semi-standard Young tableau as a corollary to Theorem {\ref{atoms}}.  Begin with a semi-standard Young tableau $T$ and map $T$ to the semi-skyline augmented filling $\Psi(T)$ whose weight is equal to that of $T$. Let the shape of $\Psi(T)$ be given by the composition $\gamma$.  Then the right key of $T$ is the unique key with weight $\gamma$.

The Demazure character $\kappa_{\omega(\lambda)}$ corresponding to a partition $\lambda$ and permutation $\omega$ can be described combinatorially as the sum of the weights of all SSYT whose right key is less than or equal to $K(\omega, \lambda)$.  (In this paper we use the notation $\kappa_{\omega(\lambda)}$ as in \cite{RS95} for ease of notation and to emphasize the fact that the Demazure characters we are working with coincide with the key polynomials described by Reiner and Shimozono.  The notation $D_{\omega}(e^{\lambda})$ typically refers to the Demazure character corresponding to a highest-weight module of weight $\omega \lambda$ over an arbitrary symmetrizable Kac-Moody Lie algebra~\cite{Kashi}.)

Demazure characters can be computed by summing over Demazure atoms.  That is, $$\kappa_{\omega(\lambda)}=\sum_{\tau \le \omega} \mathfrak{U}(\tau, \lambda),$$ where the ordering on the permutations is the Bruhat order.  A permuted-basement semi-skyline augmented filling is defined by rules similar to those which describe an ordinary semi-skyline augmented filling.  Permuting the basements of semi-skyline augmented fillings provides an alternate method for computing Demazure characters combinatorially.

\begin{theorem}{\label{keysum}}
The Demazure character $\kappa_{\omega(\lambda)}$ is equal to the sum of the weights of all permuted-basement semi-skyline augmented fillings of shape $\lambda$ with basement $\omega$.
\end{theorem}

A similar connection exists between nonsymmetric Macdonald polynomials specialized to $q=t= \infty$ and Demazure characters of the corresponding affine Kac-Moody algebra~\cite{I1}.  This correspondence and its proof provide a representation-theoretic perspective on the role of nonsymmetric Macdonald polynomials in the study of affine Lie algebras.

\section{Demazure characters}

Let $\mathfrak{g}=\mathfrak{g}l_n(\mathbb{C})$ be the general linear Lie algebra and let $\Phi$ be the corresponding root system whose highest weights are partitions.  If $\mathfrak{n}$ is the subalgebra of $\mathfrak{g}$ with basis $X_{\alpha} \; \; (\alpha \in \Phi^+)$, then $U(\mathfrak{n})$ is the universal enveloping algebra.  Let $V(\lambda)$ be the irreducible highest-weight module of weight $\lambda$.  Given a permutation $\omega$, let $u_{\omega \lambda}$ be the extremal vector of weight $\omega \lambda$ of $V(\lambda)$.  Then the character $\kappa_{\omega(\lambda)}$ of $U(\mathfrak{n})u_{\omega \lambda}$ is given by the Demazure character formula.  In this section we describe the explicit formula using Demazure operators. 

\subsection{The Demazure operator}

Let $P$ be the polynomial ring $\mathbb{Z}[x_1,x_2, \hdots]$ and let $S_{\infty}$ be the permutation group of the positive integers.  This group acts on $P$ by permuting the indices of the variables.  If $s_i$ is the elementary transposition $(i,i+1)$, define the linear operators $\partial_i$ and $\pi_i$ as in \cite{RS95} by 
\begin{equation}{\label{key_operators}}
\partial_i = \frac{1-s_i}{x_i-x_{i+1}}, \\
\; \; \; \; \; \; \; \; \; \; \pi_i = \partial_i x_i.
\end{equation}

Given $\omega \in S_{\infty}$, let $\omega=s_{i_1} s_{i_{2}} \hdots s_{i_k}$ be a decomposition of $\omega$ into elementary transpositions.  When the number $k$ of transpositions in such a product is minimized, the word $i_1 i_2 \hdots i_k$ is called a {\it reduced word} for $\omega$.  The operator $\pi_{\omega}=\pi_{i_{1}} \pi_{i_{2}} \hdots \pi_{i_k}$ is obtained by applying the product of the operators $\pi_{i_j}$, where $i_1 i_2 \hdots i_k$ is a reduced word for $\omega$.  This operator is the Demazure operator~\cite{D1},~\cite{Jo85} for the general linear Lie algebra $\mathfrak{gl}_n(\mathbb{C})$.  One obtains the {\it Demazure character} corresponding to a partition $\lambda$ and a permutation $\omega$ by applying the operator $\pi_{\omega}$ to the dominant monomial $x^{\lambda}=\prod_i x_i^{\lambda_i}$.  For example, if $\lambda=(2,1)$ and $\omega=(1,2,3)$, then the corresponding Demazure character is $\pi_{1} \pi_{2} (x_1^2 x_2)=x_1^2 x_2 + x_1^2 x_3 + x_1 x_2^2 + x_1 x_2 x_3 + x_2^2 x_3$.

\subsection{An equivalent definition~\cite{RS95}}

A {\it key} is a semi-standard Young tableau such that the set of entries in the $(j+1)^{th}$ column form a subset of the set of entries in the $j^{th}$ column, for all $j \ge 1$.  A bijection exists between weak compositions and keys given by $\gamma = (\gamma_1, \gamma_2, \hdots) \mapsto \mathit{key}(\gamma)$, where $\mathit{key}(\gamma)$ is the key such that for all $j$, the first $\gamma_j$ columns contain the letter $j$.  To invert this map, send the key $T$ to the composition describing the content of $T$.  Figure {\ref{wtf}} depicts $\mathit{key}(2,1,1,4,0,3)$, written in French notation.

\begin{figure}
$$T=\tableau{6 \\ 4 \\ 3 & 6 \\ 2 & 4 & 6 \\ 1 & 1 & 4 & 4}$$
\caption{$T=\key(2,1,1,4,0,3)$}
{\label{wtf}}
\end{figure}

Let $\col(T)$ be the word obtained from an SSYT $T$ by reading the column entries of $T$ from top to bottom, left to right.  If a word $w$ is {\it Knuth equivalent} \cite{Knu70} to $\col(T)$, write $w \sim T$.  There exists a unique word $v$ in each Knuth equivalence class such that $v = \mathit{col}(T)$ for some semi-standard Young tableau $T$.  

The {\it column form} of a word $w$, denoted $\colform(w)$, is the composition consisting of the lengths of the strictly decreasing subwords of $w$.  Let $w$ be an arbitrary word such that $w \sim T$ for $T$ of shape $\lambda$.  The word $w$ is said to be {\it column-frank} if $\mathit{colform}(w)$ is a rearrangement of the nonzero parts of $\lambda'$, where $\lambda'$ is the conjugate shape of the partition $\lambda$ obtained by reflecting the Ferrers diagram of $\lambda$ across the line $x=y$.  (In Figure {\ref{frank}}, $v$ is not column-frank but $w$ is.)

\begin{figure}[b]
{\label{frank}}
\mbox{\subfigure{
\begin{picture}(100,100)
\put(0,80){$$S = \tableau{5 \\ 3 \\ 2 & 4 \\ 1 & 2}$$}
\put(0,20){$S \sim v=3 \; 5 \; 4 \; 2 \; 2 \; 1$}
\put(0,0){$\colform(v)=(1,3,2)$}
\end{picture}}
\subfigure{
\begin{picture}(100,100)
\put(0,80){$$T = \tableau{4 \\ 3 \\ 2 & 5 \\ 1 & 2}$$}
\put(0,20){$T \sim w= 4 \; 2 \; 5 \; 3 \; 2 \; 1$}
\put(0,0){$\colform(w)=(2,4)$}
\end{picture}
}
}
\caption{Here $v$ is not column-frank but $w$ is.}
\end{figure}

Let $T$ be a semi-standard Young tableau of shape $\lambda$.  The {\it right key} of $T$, denoted $K_{+}(T)$, is defined in \cite{RS95} to be the unique key of shape $\lambda$ whose $j^{th}$ column is given by the last column of any column-frank word $v$ such that $v \sim T$ and $\mathit{colform}(v)$ is of the form $(\hdots, \lambda_j')$.  For example, $S= 5 \; 3 \; 2 \; 1 \cdot 4 \; 2$ in Figure {\ref{frank}} has right key $K_+(T)= 5 \; 4 \; 2 \; 1 \cdot 4 \; 2$.

Given an arbitrary partition $\lambda$ and permutation $\omega$ (written in one-line notation), there exists an associated key ${K}(\omega, \lambda)$ defined as follows.  Consider the subword consisting of the first $\lambda_1$ letters of $\omega$ and reorder the letters in decreasing order.  This is the first column of ${K}(\omega, \lambda)$.  The second column of ${K}(\omega, \lambda)$ contains the first $\lambda_2$ letters of $\omega$ in decreasing order.  Continuing this way, one derives the word $\mathit{col}({K}(\omega, \lambda))$ \cite{LS90}.  For example, $\omega=241635$ and $\lambda=(4,2,2,1)$ give the key $K(\omega,\lambda)=6 \; 4 \; 2 \; 1 \cdot 4 \; 2 \cdot 4 \; 2 \cdot 2$.

Define a partial order on the set of all semi-standard Young tableaux of shape $\lambda$ by setting $T \le S$ if and only if the entry in the $i^{th}$ row and $j^{th}$ column of $T$ is less than or equal to the corresponding entry in $S$ for all $i$ and $j$.  The key polynomial $\kappa_{\omega(\lambda)}$ is defined \cite{RS95} as the sum of the weights of all SSYT having right key less than or equal to ${K}(\omega, \lambda)$.  This polynomial is precisely the type $A$ Demazure character $\pi_{\omega} (x^{\lambda})$ \cite{LS90}, so we use these terms interchangeably.

\subsection{Intersections of key polynomials}

Notice that for a fixed partition $\lambda$, the sets of semi-standard Young tableaux contributing weights to the polynomials $\pi_{\omega} (x^{\lambda})$ intersect nontrivially.  For example, $$\pi_{(1,2,3)}(x^{(2,1)})=\pi_{1} \pi_{2} (x_1^2 x_2)=x_1^2 x_2 + x_1^2 x_3 + x_1 x_2^2 + x_1 x_2 x_3 + x_2^2 x_3,$$ $$\pi_{2}(x^{(2,1)})=\pi_2(x_1^2 x_2) = x_1^2 x_2 + x_1^2 x_3,$$ $$\pi_{1}(x^{(2,1)})=\pi_1(x_1^2 x_2) = x_1^2 x_2 + x_1 x_2^2$$  Here the monomial $x_1^2 x_3$ appears in both $\pi_{1} \pi_{2} (x_1^2 x_2)$ and $\pi_{2}(x^{(2,1)})$, but the SSYT $T$ with column word $\col(T)=3 \; 1 \cdot 1$ is the only SSYT of shape $\lambda=(2,1)$ and weight $x_1^2x_3$.

In fact, the definition of the operator $\pi_{\omega}$ implies that if $\omega$ is longer than $\sigma$ and $\omega=s_i \sigma $, then the semi-standard Young tableaux appearing as monomials in $\pi_{\sigma} (x^{\lambda})$ are a subset of those appearing in $\pi_{\omega} (x^{\lambda}) = \pi_i \pi_{\sigma} (x^{\lambda})$.  Therefore it makes sense to consider the intersections and complements of Demazure characters.

Let $\omega$ be a permutation of length $k$.  Consider all permutations $\sigma$ less than $\omega$ in the Bruhat order.  We study the subset of monomials in $\pi_{\omega}(x^{\lambda})$ which do not appear in $\pi_{\sigma} (x^{\lambda})$ for any such $\sigma$.  The sum of these monomials is the polynomial obtained by replacing the operators $\pi_i$ by the operators $\overline{\pi_i}=\pi_i-1$ in the formula $\pi_{\omega}(x^{\lambda})$.  The operator $\overline{\pi}_i = \overline{\pi}_{s_i}$ is therefore defined by

\begin{displaymath}
f \longrightarrow (s_i (f) - f)/(1-x_i/x_{i+1}) =  \overline{\pi}_i (f),
\end{displaymath}
and, given any reduced word $s_{i_1} s_{i_2} \hdots s_{i_k}$ for $\omega$, define $\overline{\pi}_{\omega} (f) = \overline{\pi}_{i_1} (f)  \overline{\pi}_{i_2} (f) \hdots  \overline{\pi}_{i_k} (f)$.  For example, if $f=x_1^2x_2x_3$, then $\overline{\pi}_1 f = \frac{(x_1x_2^2x_3- x_1^2x_2x_3)}{(1-x_1/x_2)}=x_1x_2^2x_3$.

Lascoux and Sch\"{u}tzenberger \cite{LS90} call these polynomials the {\it standard bases} and prove that the standard basis $\mathfrak{U}(\omega,\lambda)$ equals the sum of the weights of all SSYT having right key equal to ${K}(\omega, \lambda)$.  We retain the notation $\mathfrak{U}(\omega,\lambda)$ but call the polynomials {\it Demazure atoms} to avoid confusion with various objects referred to as standard bases.

The operators $\overline{\pi}_i$ satisfy the Coxeter relations $\overline{\pi}_i\overline{\pi}_{i+1}\overline{\pi}_i=\overline{\pi}_{i+1}\overline{\pi}_i\overline{\pi}_{i+1}$ and $\overline{\pi}_i\overline{\pi}_j=\overline{\pi}_j\overline{\pi}_i$ for $\| j-i \| > 1$~\cite{LS90}.  Lift the operator $\overline{\pi}_i$ to an operator $\theta_i$ on the free algebra by the following process.  Given $i$ and a word $w$ in the commutative alphabet $X =( x_1, x_2, \hdots$), let $m_j$ be the number of occurrences of the letter $x_j$ in $w$, for each $j$.  Let $k=m_i-m_{i+1}$.  If $k \ge 0$, then $w$ and $w^{s_i}$ differ by the exchange of a subword $x_i^k$ with the subword $x_{i+1}^k$.  The analogous statement is true for $k < 0$.  When $k \ge 0$, define $w\theta_i$ to be the sum of all words in which the subword $x_i^k$ of $w$ has been changed respectively into $x_i^{k-1}x_{i+1}, \; \; x_i^{k-2}x_{i+1}^2, \; \; \hdots , \; \; x_{i+1}^k$.  For example, if $w=x_1^3 x_2^4 x_3  x_5 x_7^3$, then $w \theta_2 = x_1^3 x_2^3 x_3^2  x_5 x_7^3 + x_1^3 x_2^2 x_3^3 x_5 x_7^3 + x_1^3 x_2 x_3^4  x_5 x_7^3$.

Every partition $\lambda=(\lambda_1,\lambda_2,\hdots)$ has a corresponding {\it dominant monomial}, $$x^{\lambda} =\prod_i x_i^{\lambda_i},$$ which equals the weight of the {\it Yamanouchi tableau} of shape $\lambda$.  (The Yamanouchi tableau is the SSYT such that, for each $i$, the entries in the $i^{th}$ row are all equal to $i$.)

\begin{theorem}{(Lascoux-Sch\"{u}tzenberger \cite{LS90})}
{\label{defthm}}
Let $x^{\lambda}$ be the dominant monomial corresponding to $\lambda$ and let $s_{i_1} s_{i_2} \hdots s_{i_k}$ be any reduced decomposition of a permutation $\pi$.  Then $\mathfrak{U}(\omega,\lambda)= x^{\lambda} \theta_{i_1} \theta_{i_2} \hdots \theta_{i_k} $.
\end{theorem}

Theorem {\ref{defthm}} provides an inductive method for constructing the Demazure atom $\mathfrak{U}(\omega,\lambda)$.  Begin with $\mathfrak{U}(id,\lambda)=x^{\lambda}$ and apply $\theta_{i_1}$ to determine $\mathfrak{U}(s_{i_1},\lambda)$.  Then apply $\theta_{i_2}$ to $\mathfrak{U}(s_{i_1},\lambda)$ to determine $\mathfrak{U}(s_{i_1} s_{i_2},\lambda)$.  Continue this process until the desired standard basis $\mathfrak{U}(\omega,\lambda)$ is obtained.

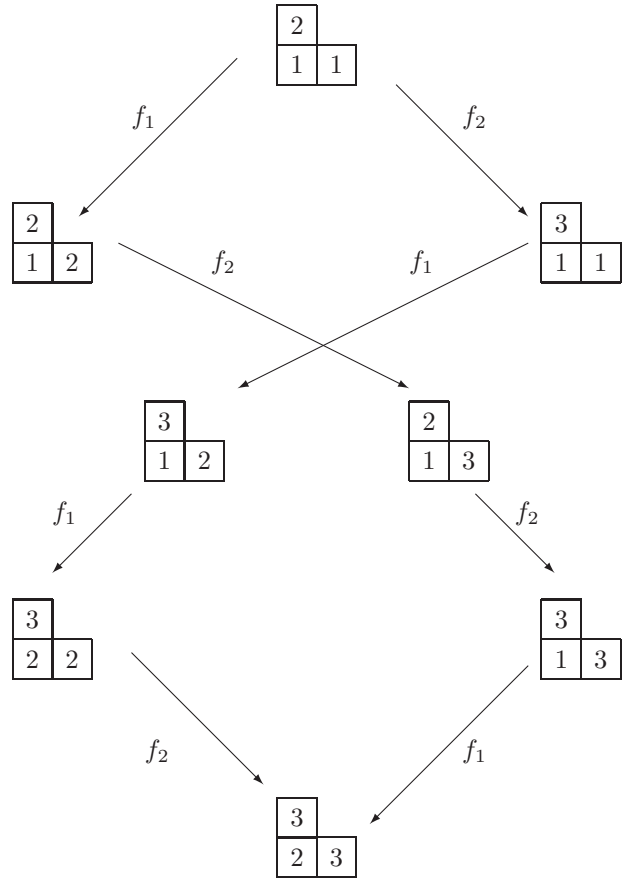
\begin{figure}[!]
{\label{base}}
\begin{center}
\begin{picture}(200,320)
\put(100,305){$$\tableau{2 \\ 1 & 1}$$}
\put(200,230){$$\tableau{3 \\ 1 & 1}$$}
\put(0,230){$$\tableau{2 \\ 1 & 2}$$}
\put(50,155){$$\tableau{3 \\ 1 & 2}$$}
\put(150,155){$$\tableau{2 \\ 1 & 3}$$}
\put(0,80){$$\tableau{3 \\ 2 & 2}$$}
\put(200,80){$$\tableau{3 \\ 1 & 3}$$}
\put(100,5){$$\tableau{3 \\ 2 & 3}$$}
\put(145,290){\vector(1,-1){50}}
\put(45,75){\vector(1,-1){50}}
\put(175,135){\vector(1,-1){30}}
\put(85,300){\vector(-1,-1){60}}
\put(195,70){\vector(-1,-1){60}}
\put(45,135){\vector(-1,-1){30}}
\put(40,230){\vector(2,-1){110}}
\put(195,230){\vector(-2,-1){110}}
\put(170,275){$f_2$}
\put(45,275){$f_1$}
\put(150,220){$f_1$}
\put(75,220){$f_2$}
\put(15,125){$f_1$}
\put(190,125){$f_2$}
\put(50,35){$f_2$}
\put(170,35){$f_1$}
\end{picture}
\caption{The crystal graph for $\lambda=(2,1)$.}
\end{center}
\end{figure}

Lascoux and Sch\"{u}tzenberger further break down this procedure to produce a crystal graph structure~\cite{LS90}.  (Throughout this paper, our crystallographic notation will follow the notation appearing in \cite{K1}.)  To describe the operator $f_i$ needed for this procedure, let $\col(T)$ be the column word corresponding to the semi-standard young tableau $T$.  Change all occurrences of $i$ in $\col(T)$ to right parentheses and all occurrences of $i+1$ in $\col(T)$ to left parentheses.  Ignore all other entries in $\col(T)$ and match the parentheses in the usual manner.  If there are no unmatched right parentheses, then $f_i(\col(T))=\col(T)$.  Otherwise replace the rightmost unmatched right parenthesis by a left parenthesis and convert the parentheses back to occurrences of $i$ and $i+1$.  The resulting word is $f_i(\col(T))$.  Figure {\ref{base}} depicts the crystal graph corresponding to the partition $(2,1)$.

The Demazure character corresponding to $\omega=s_{i_1} s_{i_2} \hdots s_{i_k}$ is obtained from this procedure by applying the appropriate $f_i$ operators.  To see this, begin with the SSYT of highest weight, which corresponds to the monomial $x^{\lambda}$.  Apply $f_{i_k}^{m_{i_k}}$, where $m_{i_k}$ is the number of unmatched right parantheses.  Add the resulting monomials to the initial monomial to obtain the Demazure character $\kappa_{s_{i_k}(\lambda)}$.  Next apply $f_{i_{k-1}}^{m_{i_{k-1}}}$ to the monomials in $\kappa_{s_{i_k}(\lambda)}$ and collect these monomials together with the monomials of $\kappa_{s_{i_k}}$ to obtain $\kappa_{s_{i_{k-1}} s_{i_k}(\lambda)}$.  Continue this procedure to obtain $\kappa_{\omega(\lambda)}$.

\section{Combinatorial description of $E_{\gamma}(X;0,0)$}

The polynomials $E_{\gamma}(X;0,0)$ are obtained from the nonsymmetric Macdonald polynomials by letting $q$ and $t$ approach $0$.  The combinatorial formula for nonsymmetric Macdonald polynomials provided by Haglund, Haiman, and Loehr \cite{HHL05c} can be specialized in this manner to obtain a combinatorial formula for $E_{\gamma}(X;0,0)$.  Several definitions are needed to describe this formula.

Let $\gamma=(\gamma_1,\gamma_2,\hdots)$ be a weak composition of $n$.    The {\it column diagram of $\gamma$} is a figure $dg'(\gamma)$ consisting of $n$ {\it
cells} arranged into columns, as in \cite{HHL05c}.   The $i^{th}$ column contains $\gamma_i$ cells, and the number of cells in a column is called the {\it
height} of that column.  A cell $a$ in a column diagram is denoted $a=(i,j)$, where $i$ is the row and $j$ is the column of the cell containing $a$.  

For example, the following depicts the column diagram of $\gamma=(0,2,0,3,1,2,0,0,1)$.
$$dg'(\gamma)=\tableau{ & & & {} \\ & {} & & {} & & {} \\ & {} & & {} & {} & {} & & & {}}$$

The {\it augmented diagram} of $\gamma$, defined by $\widehat{dg}(\gamma)=dg'(\gamma) \cup \{ (0,i) : 1 \le i \le m \}$ (where $m$ is the number of parts of $\gamma$), is the column diagram with $m$ extra cells adjoined in row $0$.  In this paper the adjoined row, called the {\it basement}, always contains the numbers $1$ through $m$ in strictly increasing order. 

The augmented diagram for $\gamma=(0,2,0,3,1,2,0,0,1)$ is depicted below.
$$\widehat{dg}(\gamma)=\tableau{ & & & {} \\ & {} & & {} & & {} \\ & {} & & {} & {} & {} & & & {} \\ 1 & 2 & 3 & 4 & 5 & 6 & 7 & 8 & 9}$$

An {\it augmented filling,} $\sigma$, of an augmented diagram $\widehat{dg}(\gamma)$
is a function $\sigma: \widehat{dg}(\gamma) \rightarrow \mathbb{Z}_+$, which we
picture as an assignment of positive integer entries to the cells of
$\gamma$.  Let $\sigma(k)$ denote
the entry in the $k^{th}$ cell of the augmented diagram
encountered when $\widehat{dg}(\gamma)$ is read across rows from left to right, beginning at
the highest row and working downward.  This ordering of the cells is called the {\it reading order}.  (A cell $a=(i,j)$ is greater than a cell $b=(i',j')$ in the reading order if either $i'>i$ or $i'=i$ and $j'<j$.)  The reading word $\mathit{read}(\sigma)$ is obtained by recording the entries in this reading order.  The {\it content} of a filling $\sigma$ is the multiset of entries which appear in the filling.  The cells $a_1=(i_1,j_1)$ and $a_2=(i_2,j_2)$ of $\gamma$ are said to be {\it attacking} if any of the following three conditions are true:
\begin{itemize}
\item $i_1 = i_2$
\item $i_1 - i_2 = 1$ and $j_2 < j_1$
\item $i_2 - i_1 = 1$ and $j_1 < j_2$.
\end{itemize}
A filling is said to be {\it non-attacking} if for every pair of attacking cells $\{ a_1, a_2 \}$, we have $\sigma(a_1) \not= \sigma(a_2)$.  The fillings utilized in the combinatorial description of Demazure atoms are non-attacking fillings with additional row and column restrictions.

The following triples of cells are introduced to provide restrictions on the row entries of a filling.  Note that the triple types are not related to symmetry types.  Type $A$ and type $B$ merely refer to the positions of the cells in the diagram.

Let $a_1=(i_1,j_1), a_2 = (i_2,j_2),$ and $a_3 =(i_3,j_3)$ be three cells in $\widehat{dg}(\gamma)$ such that column $j_1$ is taller than or equal in height to column $j_2$.  If $i_1=i_2$, $i_1-i_3=1$, and $j_1=j_3$, then $a_1,a_2,$ and $a_3$ are said to form a {\it type $A$ triple}, as depicted below.
$$\tableau{a_1 \\ a_3} \cdots \tableau{a_2}$$

Define for $x,y \in \mathbb{Z}_+$
\[I(x,y) = \left\{ \begin{array}{ll}
1  & \mbox{if $x>y$} \\
0  & \mbox{if $x \le y$}
\end{array}
\right. . \] 
Let $\sigma$ be an augmented filling and let $\{ \sigma(a_1),\sigma(a_2),\sigma(a_3) \}$ be the entries of
$\sigma$ in the cells $\{a_1,a_2,a_3\}$, respectively, of a type $A$ triple.  The triple $\{a_1,a_2,a_3\}$ is called a {\it type $A$ inversion triple} if and only
if $I(\sigma(a_1),\sigma(a_2))+I(\sigma(a_2),\sigma(a_3))-I(\sigma(a_1),\sigma(a_3))=1$.

Consider the following ordering of the cells $a_1,a_2,a_3$ obtained from their entries.  Let $a_i < a_j$ if either $\sigma(a_i) < \sigma(a_j)$ or $\sigma(a_i) = \sigma(a_j)$ and $a_i$ comes before $a_j$ in reading order.  If this ordering produces a counter-clockwise orientation of the cells $a_1,a_2,a_3$ when read from smallest to largest, then the cells form a type $A$ inversion triple.  This definition is equivalent to that given by the function $I(x,y)$.

Similarly, consider three cells $\{a_1=(i_1,j_1,a_2=(i_2,j_2),a_3=(i_3,j_3) \} \in \lambda$ such that column $j_2$ is strictly taller then column $j_1$.  The cells $\{a_1,a_2,a_3 \}$ are said to form a {\it type $B$ triple} if $i_1=i_2$, $j_2=j_3$, and $i_3-i_2=1$, as shown below.

$$\tableau{ \\ a_1} \hdots  \tableau{a_3 \\ a_2}$$

Let $\sigma$ be an augmented filling and let $\{ \sigma(a_1), \sigma(a_2), \sigma(a_3) \}$ be the entries of
$\sigma$ in the cells $\{a_1,a_2,a_3 \}$ of a type $B$ triple.  The triple $\{ a_1,a_2,a_3 \}$ is called a {\it type $B$ inversion triple} if and
only if $I(\sigma(a_3),\sigma(a_1))+I(\sigma(a_1),\sigma(a_2))-I(\sigma(a_3),\sigma(a_2))=1$.  

As for type $A$ inversion triples, there is an equivalent definition for type $B$ inversion triples.  Again let $a_i < a_j$ if either $\sigma(a_i) < \sigma(a_j)$ or $\sigma(a_i) = \sigma(a_j)$ and $a_i$ comes before $a_j$ in reading order.  The three cells form a type $B$ inversion triple if the ordering of the cells, when read from smallest to largest, produces a clockwise orientation.

Define a {\it semi-skyline augmented filling} of an augmented diagram $\widehat{dg}(\gamma)$ to be an augmented filling $F$ such that the entries in each column (read top to bottom) are weakly increasing and every type $A$ or type $B$ triple of cells is an inversion triple.  Corollary 2.4 of \cite{M06} states that these conditions are enough to guarantee that the filling is non-attacking.  Specializing the combinatorial formula for the nonsymmetric Macdonald polynomials $E_{\gamma}(x;q,t)$ given in \cite{HHL05c} implies that $$E_{\gamma}(x;0,0)=\sum_{F \in SSAF(\widehat{dg}(\gamma))} x^F,$$ where $SSAF(\widehat{dg}(\gamma))$ is the set of all semi-skyline augmented fillings of shape $\gamma$.

\section{Proof of Theorem \ref{atoms}}

The set of Demazure atoms for the partition $\lambda$ can be considered as a decomposition of the Schur function $s_{\lambda}$.  For any partition $\lambda$ of $n$, it is known \cite{LS90} that $$\sum_{\omega \in S_n} \mathfrak{U}(\omega,\lambda)=s_{\lambda}.$$  The functions $E_{\gamma}(X;0,0)$ are also a decomposition of the Schur functions \cite{M06}, so it is natural to determine their relationship to the Demazure atoms.  Theorem {\ref{atoms}} states that $\mathfrak{U}(\omega,\lambda) = NS_{\omega(\lambda)}$, where $\omega(\lambda)$ denotes the action of $\omega$ on the parts of $\lambda$ when $\lambda$ is considered as a partition of $n$ into $n$ non-negative parts. 

\subsection{Several Useful Lemmas}

Section 3.1 of \cite{M06} provides a bijection $\rho$ between row-strict plane partitions  and semi-skyline augmented fillings which preserves the entries in each row.  The map can be considered as a map from a collection $\{ R_i \}$ of sets of row entries to an SSAF.  Insert the rows from lowest to highest.  Assume that the lowest $j$ rows and the largest $k$ entries of row $j+1$ have been inserted.  Consider $\alpha_{k+1}$, the $(k+1)^{th}$ largest entry in row $j+1$.  Place $\alpha_{k+1}$ on top of the leftmost entry, $\beta$, of row $j$ such that the cell on top of $\beta$ is empty and $\beta \ge \alpha_{k+1}$.  Continue in this manner until all the row entries have been placed into the diagram.  The result is the unique SSAF with row entries $\{ R_i \}$.

A different bijection, $\Psi$, is described in \cite{M06} to map directly between SSYTs and SSAFs.  Begin with a semi-standard Young tableau $T$ of shape $\lambda$ and insert its entries into an empty SSAF, using the following insertion procedure.  When inserting an entry $\alpha_1$ into an SSAF $F$, find the first entry $\alpha_2$ of $F$ in reading order which is greater than or equal to $\alpha_1$.  If there is no entry on top of $\alpha_2$, place $\alpha_1$ on top of $\alpha_2$ and the insertion is complete.  If the entry directly above $\alpha_2$ is greater than $\alpha_1$, continue to the next entry in reading order which is greater than or equal to $\alpha_1$ and repeat.  If the entry, $\alpha_3$, directly above $\alpha_2$ is less than $\alpha_1$, replace it with $\alpha_1$ and find the next entry in reading order which is greater than or equal to $\alpha_3$.  Repeat this procedure until the insertion is complete.  Applying this insertion to the columns of $T$, beginning with the smallest entry in the rightmost column and moving through the column entries from smallest to largest, rightmost column to leftmost column, produces a semi-skyline augmented filling whose shape is a rearrangement of $\lambda$.  This map is a weight-preserving, shape-rearranging bijection between semi-standard Young tableaux and semi-skyline augmented fillings \cite{M06}.

\begin{proposition}
{\label{diagram}}
There exists a map $\Theta_i :$SSAF$\longrightarrow$SSAF such that the following diagram commutes for all SSYT $T$.

\begin{center}
\begin{picture}(100,40)
\put(0,40){\xymatrix{T \ar[r]^{\mbox{$f_i$}} \ar[d]_{\mbox{$\Psi$}} & T' \ar[d]^{\mbox{$\Psi$}} \\ F \ar[r]^{\mbox{$\Theta_i$}} & F'}}
\end{picture}
\end{center}
(Here $f_i$ is the crystal operator described in Section 2.)

\end{proposition}

\begin{proof}
Let $F$ be an arbitrary semi-skyline augmented filling and let $\mathit{read}(F)$ be the reading word obtained by reading $F$ left to right, top to bottom, as described in Section 3.  First match any pair $i$ and $i+1$ which occur in the same row of $F$ and remove these entries from the reading word of $F$.  Next apply the parenthetical matching procedure of \cite{LS90} described in Section 2.3 to the reading word to determine which of the remaining occurrences of $i$ and $i+1$ are unmatched.    In other words, replace each $i+1$ by a left (open) parenthesis and each $i$ by a right (closed) parenthesis and match left and right parenthesis.  

Pick the rightmost unmatched $i$.  Convert it to an $i+1$.  (If there is no unmatched $i$, then $\Theta_i(F)=F$.)  The result is a collection of row entries which differ from those of $\mathit{read}(F)$ in precisely one entry.  Use the procedure $\rho$ described above to map this collection of rows to a unique SSAF.  (This map is well-defined for our collection of row entries because the row directly below the rightmost unmatched $i$ either does not contain the entry $i$ or contains both the entry $i$ and the entry $i+1$.)  The resulting SSAF is $\Theta_i(F)=F'$.  We must show that $\Theta_i(\Psi(T))=\Psi(\tilde{f_i}(T))$.

Recall that the map $\rho$ is a bijection between semi-skyline augmented fillings and row-strict plane partitions which preserves the entries in each row.  The inverse of $\rho$ sends the entries in the $r^{th}$ row of an SSAF $F$ to the $r^{th}$ row of a row-strict plane partition in decreasing order.  The reading word for a row-strict plane partition is given by reading the rows from left to right, top to bottom.  The matching procedure on this word and therefore the operator $f_i$ are the same as those applied to the column word of an SSYT.  This means that the resulting crystal graph is the image of the ordinary crystal graph under the weight-preserving bijection between SSYT of a fixed shape and row-strict plane partitions of the transposed shape.  Therefore it is enough to show that the following diagram commutes for all row-strict plane partitions $P$.

\begin{center}
\begin{picture}(100,50)
\put(0,40){\xymatrix{P \ar[r]^{\mbox{$f_i$}} \ar[d]_{\mbox{$\rho$}} & P \ar[d]^{\mbox{$\rho$}} \\ F \ar[r]^{\mbox{$\Theta_i$}} & F'}}
\end{picture}
\end{center}

We claim that if the the rightmost unmatched $i$ in $F$ appears in the $r^{th}$ row of $F$, then the rightmost unmatched $i$ in $\rho^{-1}(F)$ appears in the $r^{th}$ row of $\rho^{-1}(F)$.  Moreover, if there is no unmatched $i$ in $F$, then there is no unmatched $i$ in $\rho(F)$.  

Notice that the first step in the procedure for matching entries in $F$ is to match each pair of entries $\{ i,i+1 \}$ appearing in the same row of $F$.  These entries will also appear in the same row of $\rho(F)$.  Since the row entries in $\rho(F)$ appear in strictly decreasing order, the $i+1$ appears first in the reading word.  Therefore this pair of entries will be matched in $\rho(F)$.  Once these entries are matched, the remaining occurrences of $i$ and $i+1$ appear in the same order in the reading word for $F$ as in the reading word for $\rho(F)$.  Therefore the rightmost unmatched $i$ appears in the same row of $F$ as in $\rho(F)$, and if each $i$ is matched in $F$ then each $i$ is matched in $\rho(F)$.

To see that the proposition follows from this claim, first consider the situation in which there is no unmatched $i$ in $F$.  The claim implies that there is no unmatched $i$ in $\rho(F)$.  Then $\Theta_i(F)=F$ and $f_i(T)=T$.  The diagram commutes since $F=\rho(T)$.  Next let $i_r$ be the rightmost unmatched $i$ in $F$, appearing in the $r^{th}$ row of $F$.  The $r^{th}$ row of $\Theta_i(F)$ is the only row whose entries are different from $F$.  The $i$ in this row was changed to an $i+1$.  Similarly, the $r^{th}$ row of $f_i(\rho(F))$ is the only row which whose entries are different from $\rho(F)$, and the difference is an $i$ replaced by an $i+1$.  Therefore $f_i(\rho(F))$ is the image of $\Theta_i(F)$ under $\rho$ and the diagram commutes.
\end{proof}

We need one additional Lemma to prove Theorem \ref{atoms}.

\begin{lemma}
{\label{switch}}
Let $F \in SSAF(\gamma)$.  Then either $\Theta_i(F) \in SSAF(\gamma)$ or $\Theta_i(F) \in SSAF(s_i \gamma)$.
\end{lemma}

\begin{proof}
Assume that $F \in SSAF(\gamma)$.  If $\Theta_i(F)=F$, then $\Theta_i(F) \in SSAF(\gamma)$.  We must prove that when an unmatched $i$ is sent to $i+1$, the resulting semi-skyline augmented filling is either in $SSAF(\gamma)$ or in $SSAF(s_i \gamma)$.  Let $i_r$ denote the rightmost unmatched $i$ in $F$, where $r$ is the row in which $i_r$ appears.  Similarly, let $(i+1)_r$ denote the $i+1$ which replaces $i_r$ in $\Theta_i(F)$.

If $i_r$ appears in $F$ immediately above an entry greater than $i$, then $(i+1)_r$ is mapped to the same position by $\rho$.  In this case the remaining entries of the $r^{th}$ row are mapped to the same positions in $\Theta_i(F)$ as in $F$.  If the $(r+1)^{th}$ row does not contain an $i+1$, its entries are sent to the same positions as in $F$ and therefore the shape of $\Theta_i(F)$ is equal to the shape of $F$.  If there exists an $i+1=(i+1)_{r+1}$ in the $(r+1)^{th}$ row of $F$, then there must also be an $i=i_{r+1}$ since $i_r$ is unmatched.  The only situation in which the placement of this row into $\Theta_i(F)$ differs from its placement in $F$ is if $(i+1)_{r+1}$ appears to the right of $i_r$ in $F$.  In this case $(i+1)_{r+1}$ would be inserted on top of $(i+1)_r$ and $i_{r+1}$ would replace $(i+1)_{r+1}$.  The remaining entries of row $r+1$ would be inserted into the same positions in $\Theta_i(F)$ as in $F$.  The remaining rows of $\Theta_i(F)$ are inserted into the same positions as in $F$ unless row $r+2$ contains an $i$ and an $i+1$.  If row $r+2$ does contain both $i$ and $i+1$, then a similar argument shows that either $i_{r+2}$ and $(i+1)_{r+2}$ appear in the same positions as in $F$ or switch positions in $\Theta_i(F)$.  Repeating this argument for each row implies that the shape of $\Theta_i(F)$ is the same as that of $F$.  

If $i_r$ appears in $F$ immediately above another entry equal to $i$, denote this entry by $i_{r-1}$ .  Since $i_r$ is the rightmost unmatched $i$ in $read(F)$, row $r-1$ must contain an $i+1=(i+1)_{r-1}$.  This entry must appear to the right of $i_{r-1}$, for otherwise $i_r$ would appear on top of $(i+1)_{r-1}$ in $F$.    In this case the entry immediately below $i_{r-1}$ must be equal to $i=i_{r-2}$ (regardless of the column heights) in order to satisfy the inversion conditions of an SSAF.  Then there must be an $i+1$ to the right of $i_{r-2}$ in this row as well.  Applying the same arguments inductively to each row of $F$ implies that there must be an $i$ and an $i+1$ in the first row of $F$.  The semi-skyline augmented filling conditions imply that an entry $\alpha$ in the first row of $F$ must appear in the $\alpha^{th}$ column of $F$~\cite{M06}.  Therefore the entries $i$ and $i+1$ in the first row of $F$ must lie in the $i^{th}$ and $(i+1)^{th}$ columns respectively.  Thus $i_{r-1}$ appears in the $i^{th}$ column of $F$ and $(i+1)_{r-1}$ appears in the $(i+1)^{th}$ column.  The entry $(i+1)_r$ therefore passes the $i^{th}$ column and is placed into the $(i+1)^{th}$ column of $\Theta_i(F)$.

If there is no entry on top of $(i+1)_{r-1}$ in $F$, then all other entries in row $r$ of $F$ are placed in the same positions in $\Theta_i(F)$.  If an entry appears in the $(i+1)^{th}$ column of $F$, then this entry will be inserted onto the $i^{th}$ column in $\Theta_i(F)$.  In both cases, the entries in the $i^{th}$ and $(i+1)^{th}$ column are permuted and all other entries remain the same.  

Let $\beta$ be the (possibly empty) entry which lies in the $(i+1)^{th}$ column of the $r^{th}$ row of $F$ and hence the $i^{th}$ column of the $r^{th}$ row of $\Theta_i(F)$.  This entry $\beta$ must be less than $i$.  If there is no $i+1$ in row $r+1$ of $F$, then the entries on top of $\beta$ and $i+1$  in $\Theta_i(F)$ might be permuted but the other entries in row $r+1$ of $\Theta_i(F)$ retain the same positions they held in $F$.  In this case the remaining rows are the same as in $F$ up to a permutation of the entries in the $i^{th}$ and $(i+1)^{th}$ column and hence $\Theta_i(F)$ is either in $SSAF(\gamma)$ or $SSAF(s_i(\gamma))$.  If there is an $i+1$ in row $r+1$ of $F$, then there is an $i$ in row $r+1$ as well since $i_r$ is unmatched in $F$.  Then $i_{r+1}$ must lie in the $i^{th}$ column of $F$ and $(i+1)_{r+1}$ must lie to the right of $i_{r+1}$.  Therefore $(i+1)_{r+1}$ is placed on top of $(i+1)_r$ in $\Theta_i(F)$ and $i_{r+1}$ is placed in the cell which contained $(i+1)_{r+1}$ in $F$.  If an entry appears on top of $\beta$ in $F$, this same entry appears on top of $\beta$ in $\Theta_i(F)$.  All other entries in row $r+1$ of $\Theta_i(F)$ remain in the same positions as in $F$.

The entries in row $r+2$ follow a similar pattern.  If there is an $i+1=(i+1)_{r+2}$ in this row of $F$, there must also be an $i=i_{r+2}$.  Then $i_{r+2}$ occupies in $\Theta_i(F)$ the position occupied by $(i+1)_{r+2}$ in $F$.  The entry in the $(i+1)^{th}$ column of $F$ occupies the $i^{th}$ column of $\Theta_i(F)$, and $(i+1)_{r+2}$ occupies the $(i+1)^{th}$ column.  Otherwise the only entries affected are the entries in the $i^{th}$ and $(i+1)^{th}$ column, which are possibly permuted.  Eventually a row is reached which does not contain an $i+1$.  At this point the argument in the previous paragraph implies that the resulting shape of $\Theta_i(F)$ is either $\gamma$ or $s_i(\gamma)$.  
\end{proof}

Consider a semi-standard Young tableau $T$ whose weight appears in $\mathfrak{U}(\omega,\lambda)$.  We abuse notation and write $T \in \mathfrak{U}(\omega,\lambda)$.  If $s_i \omega$ is longer than $\omega$, then Lascoux and Sch\"{u}tzenberger's definition of $f_i$ implies that either $f_i(T) \in \mathfrak{U}(\omega,\lambda)$ or $f_i(T) \in \mathfrak{U}(s_i \omega,\lambda)$.  To see that the objects under consideration are the same, we must show that the operators $f_i$ act the same as the operators $\Theta_i$.

\subsection{Proof of Theorem {\ref{atoms}}}

We are now ready to prove that the Demazure atoms $\mathfrak{U}(\omega, \lambda)$ are equivalent to the polynomials $E_{\omega(\lambda)}(X;0,0)$.  We abuse notation by writing $F \in E_{\omega(\lambda)}(X;0,0)$ whenever $F$ is an SSAF of shape $\omega(\lambda)$.  This abuse is justified by the fact that the monomial $x^F$ appears in $E_{\omega(\lambda)}(X;0,0)$.

\begin{proof}
Fix a partition $\lambda$ and argue by induction on the length of the permutation $\omega$ in $\mathfrak{U}(\omega,\lambda)$.  First let $\omega$ be the identity.  Then $\mathfrak{U}(\omega,\lambda)$ is the dominant monomial.  Consider $\lambda$ as a composition of $n$ into $n$ parts by adding zeros to the right if necessary.  Each cell $a$ in the first column must have $F(a)=1$, since the columns of $F$ are weakly increasing from top to bottom and the basement entry in this column is $1$.  Each cell $b$ in the second column must have $F(b) \le 2$ since the columns are weakly increasing when read top to bottom and the basement entry is $2$.  Each cell in the second column attacks the cell immediately to its left, and therefore cannot contain the entry $1$.  This means that each cell in the second column must contain the entry $2$.  Continuing inductively, we see that for all $i$, each cell $c$ in the $i^{th}$ column must have $F(c) = i$.  To see that this is indeed an SSAF, we only need to check type $A$ triples.  But if the two cells in the left-hand column are equal and less than the cell in the right-hand column, the result is a type $A$ inversion triple.  Therefore, the $E_{\lambda}(X;0,0)=\mathfrak{U}(id,\lambda)$.

Next assume that $\mathfrak{U}(\omega,\lambda)=E_{\omega(\lambda)}(X;0,0)$, for all permutations $\omega$ of length less than or equal to $k$, for some $k \ge 0$.  (Here $\omega(\lambda)$ is the composition obtained by applying the permutation $\omega$ to the columns of $\lambda$ when $\lambda$ is considered as a composition of $n$ into $n$ parts.)  Then each permutation of length $k+1$ is obtained from a permutation of length $k$ by applying an elementary transposition $s_i$.  Let $\tau$ be an arbitrary such permutation of length $k+1$ such that $\tau=s_i \omega$ for some $\omega$ of length $k$.  The monomials in $\mathfrak{U}(\tau,\lambda)$ are obtained from the monomials of $\mathfrak{U}(\omega,\lambda)$ whose image under (possibly multiple applications of) $\tilde{\theta_i}$ is not a monomial of $\mathfrak{U}(\omega,\lambda)$.  Let $T$ be an arbitrary SSYT of $\mathfrak{U}(s_i \omega,\lambda)$ such that $T=(\tilde{\theta}_i)^m(S)$ for some SSYT $S \in \mathfrak{U}( \omega,\lambda)$ and some positive integer $m$.

Repeated application of Proposition {\ref{diagram}} implies that $\Psi((\tilde{\theta_i})^m(S))=(\Theta_i)^m (\Psi(S))$.  Since $\Psi(S) \in E_{\omega(\lambda)}(X;0,0)$ by the inductive hypothesis, Lemma {\ref{switch}} implies that either $(\Theta_i)^m(\Psi(S)) \in E_{\omega(\lambda)}(X;0,0)$ or $(\Theta_i)^m(\Psi(S)) \in E_{s_i \omega(\lambda)}(X;0,0)$.  If $(\Theta_i)^m(\Psi(S)) \in E_{\omega(\lambda)}(X;0,0)$, then $\Psi((\tilde{\theta_i})^m(S)) \in E_{\omega(\lambda)}(X;0,0)$, so $(\tilde{\theta_i})^m(S) \in \mathfrak{U}(\omega,\lambda)$ because $\mathfrak{U}(\omega,\lambda)=E_{\omega(\lambda)}(X;0,0)$ by the inductive hypothesis.  This contradicts the assumption that $(\tilde{\theta_i})^m(S) \in \mathfrak{U}(s_i \omega,\lambda)$, so $(\Theta_i)^m(\Psi(S)) \in E_{s_i \omega(\lambda)}(X;0,0)$.  Therefore, $\mathfrak{U}(s_i \omega,\lambda) \subseteq E_{s_i \omega(\lambda)}(X;0,0)$, and so $\mathfrak{U}(\tau,\lambda) \subseteq E_{\tau(\lambda)}(X;0,0)$.  

To see the reverse containment, let $F$ be a filling represented by a monomial in $E_{s_i \omega(\lambda)}(X;0,0)$.  Then $F$ is an SSAF of shape $\gamma=s_i \omega (\lambda) = \tau(\lambda)$.  Consider the smallest $j$ such that $\gamma_j < \gamma_{j+1}$.  (Such a $j$ must exist, for otherwise $F \in E_{\lambda}(X;0,0)$.)  The type $B$ inversion condition implies that all of the entries in column $j$ must be equal to $j$ and the lowest $\gamma_j+1$ entries in column $j+1$ must be equal to $j+1$.  The entry $j$ cannot appear in row $\gamma_j+1$ of $F$, because if it did it would be to the right of the $j^{th}$ column and therefore attack the $j$ in row $\gamma_j$.  This implies that the entry $j+1$ in row $\gamma_j+1$ is unmatched, and in fact it is the rightmost unmatched $j+1$ in $read(F)$, since each row below row $\gamma_j+1$ contains both a $j$ and a $j+1$.  

Consider the procedure which sends the leftmost unmatched $j+1$ in $F$ to a $j$ and then inserts the resulting row entries back into an SSAF.  This is precisely the inverse of the map $\Theta_j$.  Lemma {\ref{switch}} implies that $F'=\Theta_j^{-1}(F)$ either has shape $\tau (\lambda)$ or shape $\sigma (\lambda)$ for some permutation $\sigma$ of length $k$ such that $s_i \sigma = \tau$.  If $F'$ has shape $\sigma(\lambda)$, then the inductive hypothesis implies that $F'=\Psi(T)$ for some SSYT $T \in \mathfrak{U}(\sigma, \lambda)$.  Then $\Theta_j(F')=\Psi(\tilde{\theta_j})(T)$ by Proposition {\ref{diagram}}.  Since $(\tilde{\theta_j})(T)$ is either in $\mathfrak{U}( \sigma, \lambda)$ or $\mathfrak{U}(\tau, \lambda)$ and $F=\Psi(\tilde{\theta_j})(T)$ is not in $\mathfrak{U}(\sigma, \lambda)$, then $F= \Psi(\tilde{\theta_j})(T) \in \mathfrak{U}(\tau, \lambda)$.

Apply $\Theta_j^{-1}$ until $F^{(m)}=(\Theta_j^{-1})^{(m)}(F)$ has shape $\sigma(\lambda)$.  This occurs for some $m$ less than or equal to the number of unmatched $(j+1)'s$ in $F$.  To see this, let $m$ be the number of unmatched $(j+1)'s$ in $F$ and assume that $F^{(m-1)}=(\Theta_j^{-1})^{(m-1)}(F)$ has shape $\tau$.  Since $ \tau \not= \sigma$, column $j$ is strictly shorter than column $j+1$.  Apply $\Theta_j$ to $F^{(m-1)}$ to map the $j+1$ in row $\gamma_j+1$ to $j=j_0$.  Then $j_0$ lies in the $j^{th}$ column of $F^{(m)}$ but the other entries in this row of $F^{(m)}$ remain in the same positions as in $F^{(m-1)}$.  This arrangement of the entries in row $\gamma_j+1$ implies that column $j+1$ of $F^{(m)}$ has height $\gamma_j$, and hence the shape of $F^{(m)}$ is different from the shape of $F^{(m-1)}$.  Since the shape of $F^{(m)}$ must be equal to either $s_j \sigma$ or $\sigma$, the shape of $F^{(m)}$ must be $\sigma$.  The inductive hypothesis implies that $F^{(m)}=\Psi(T)$ for some $T \in \mathfrak{U}(\sigma, \lambda)$, since $\sigma$ has length less than that of $\tau=s_j \sigma$.  Proposition {\ref{diagram}} implies that $F=\Psi(\tilde{\theta_j})^{m}(T)$, which means that $F \in \mathfrak{U}(\tau, \lambda)$.  Therefore $E_{\tau(\lambda)}(X;0,0) \subseteq \mathfrak{U}(\tau,\lambda)$.

The above shows that $E_{\tau(\lambda)}(X;0,0) = \mathfrak{U}(\tau,\lambda)$ for an arbitrary choice of permutation $\tau$ of length $k+1$.  Therefore it is true for all permutations of length $k+1$.  Applying the principle of mathematical induction  completes the proof.
\end{proof}

Theorem {\ref{atoms}} provides a non-inductive construction of the Demazure atoms.  In particular, given a partition $\lambda \vdash n$ and a permutation $\omega \in S_n$, first consider $\lambda$ as a composition of $n$ into $n$ parts by appending zeros if necessary.  Then apply the permutation $\omega$ to the columns of $\lambda$ to obtain the shape $\omega(\lambda)$.  Finally, determine all semi-skyline augmented fillings of the shape $\omega(\lambda)$.  The monomials given by the weights of these SSAFs are the monomials appearing in the Demazure atom $\mathfrak{U}(\omega,\lambda)$.

\section{Computation of right keys}

Recall that the Demazure atom $\mathfrak{U}(\omega,\lambda)$ is equal to the sum of the weights of all SSYT with right key ${K}(\omega, \lambda)$.  Therefore all of the SSYT which map to an SSAF of shape $\omega(\lambda)$ have the same right key, ${ K}(\omega, \lambda)$.  Theorem {\ref{atoms}} provides a simple method to determine the right key of a semi-standard Young tableau. The {\it super SSAF} (denoted $\mathit{super}(\gamma)$) of a composition $\gamma$ is the SSAF of shape $\gamma$ whose $i^{th}$ column contains only the entries $i$.  The weight of this SSAF is the dominant monomial in the polynomial $\mathfrak{U}(\omega,\lambda)$ under lexicographic ordering. 

\begin{cor}{\label{key}}
Given an arbitrary SSYT $T$, let $\gamma$ be the shape of $\Psi(T)$.  Then $K_+(T)=key(\gamma)$.
\end{cor}

\begin{proof}
We must show that the map $\Psi: SSYT \rightarrow SSAF$ sends a key $T$ to $\mathit{super}(\gamma)$, where $\gamma$ is the composition $\mathit{content}(T)$.  We prove this by induction on the number of columns of $T$.  If $T$ has only one column, $C_1= \alpha_1 \; \alpha_2 \; \hdots \alpha_l$, then this column maps to a filling $F$ with one row such that the $\alpha_i^{th}$ column contains the entry $\alpha_i$ for each $i$.  This is precisely $\mathit{super}(\mathit{content}(T))$.

Next assume that $\Psi(T) = \mathit{super}(\mathit{content}(T))$ for all keys $T$ with less than or equal to $m-1$ columns.  Let $S$ be a key with $m$ columns.  After the insertion of the rightmost $m-1$ columns, the figure is $\mathit{super}(\mathit{content}(S \setminus C_1))$ by the inductive hypothesis.  We must show that the insertion of the leftmost column produces $\mathit{super}(\mathit{content}(S))$.  

Let $\mathit{super}(\mathit{content}(S \setminus C_1))$ have shape $\gamma=(\gamma_1, \gamma_2, \hdots , \gamma_n)$.  For each $i$, if $\gamma_i \not= 0$, then $i=\alpha$ is an entry in $C_1$ since $S$ is a key.  Before the insertion of $C_1$, each of the cells in the $i^{th}$ column each contain the entry $i$.  Therefore $\alpha$ cannot be bumped further in the reading order than the cell in row $\gamma_i+1$ of column $i$.   This implies that each of the non-zero entries of $\gamma$ in $C_1$ must appear at or above their respective columns, and the pigeon-hole principle therefore implies that each appears at the top of its column.  The new columns are created by the entries in $C_1$ which do not appear in any of the subsequent columns of $S$.  Therefore the result is indeed the SSAF $\mathit{super}(\mathit{content}(S))$.

Since $\mathit{super}(\omega(\lambda)) \in \mathfrak{U}(\omega, \lambda)$ and $\mathfrak{U}(\omega, \lambda)$ is a collection of all SSYT with the same right key, each SSYT which maps to a SSAF of shape $\omega(\lambda)$ has the same right key as $\mathit{super}(\omega(\lambda))$.  Therefore, if $T \in \mathfrak{U}(\omega,\lambda)$, then $K_+(T)= \mathit{key}(\omega(\lambda))$.
\end{proof}

\begin{figure}
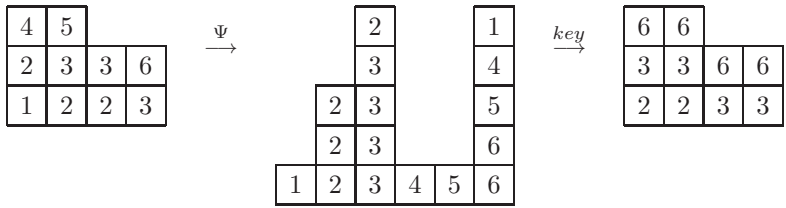

{\label{keypic2}}
$$\tableau{4 & 5 \\ 2 & 3 & 3 & 6 \\ 1 & 2 & 2 & 3} \hspace*{.2in} \substack{ \Psi \\ \longrightarrow} \hspace*{.2in} \tableau{ & & 2 & & & 1 \\ & & 3 & & & 4 \\ & 2 & 3 & & & 5 \\ & 2 & 3 & & & 6 \\ 1 & 2 & 3 & 4 & 5 & 6}  \hspace*{.2in} \substack{key \\ \longrightarrow} \hspace*{.2in} \tableau{6 & 6 \\ 3 & 3 & 6 & 6 \\ 2 & 2 & 3 & 3}$$
\caption{Computation of the key of an SSYT}
\end{figure}

Corollary {\ref{key}} provides a quick procedure for calculating the right key of any SSYT.  In particular, if $T$ is an arbitrary SSYT, then the right key of $T$ is given by $\mathit{key}(shape(\Psi(T)))$.  (See Figure {\ref{keypic2}} for an example.)  This calculation facilitates the computation of Demazure atoms and Demazure characters.  Let $\gamma$ be a composition which rearranges a partition $\lambda$, so that $\gamma = \omega(\lambda)$.  The {\it key polynomial} $\kappa_{\gamma}$ is given by the weights of all semi-standard Young tableaux $T$ of shape $\lambda$ such that $K_+(T) \le key(\gamma)$~\cite{LS90}.  Therefore $$\kappa_{\gamma}=\sum_{\alpha \le \gamma} NS_{\alpha},$$ where $\alpha \le \gamma$ if and only if $\omega_1(\alpha) = \lambda$ for some permutation $\omega_1$ such that $\omega_1 \le \omega$ in the Bruhat order.   

\section{Key polynomials and permuted basements}

Several different methods for computing key polynomials are described in~\cite{RS95}.  The notion of a semi-skyline augmented filling with a permuted basement provides an additional method which utilizes the action of the permutation group in a natural way.

Let $\gamma$ be a weak composition of $n$ into $n$ parts such that $\gamma$ is obtained from a partition $\lambda$ by applying the permutation $\omega \in S_n$.  (Here $\lambda$ is a partition of $n$ into $n$ non-negative parts.)  First construct the augmented diagram associated to $\lambda$.  Next apply the permutation $\omega$ to the entries in the basement, without permuting the columns.  Fill the cells with positive integers in such a way that the columns are weakly increasing when read top to bottom and every triple is an inversion triple.  The result, $\tilde{F}$ is called a {\it permuted basement SSAF} with shape $\lambda$ and basement $\omega$.  See Figure 6.1 for an example with the basement entries in bold.  To condense notation, we write {\it $pb(\lambda,\omega)$} to denote the set of all permuted basement SSAFs with shape $\lambda$ and basement $\omega$.  Let $\tilde{K}_{\omega, \lambda}=\sum_{\tilde{F}} x^{\tilde{F}}$ be the sum of the weights of all permuted basement SSAFs in $pb(\lambda, \omega)$.

\begin{figure}
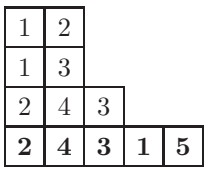

{\label{perm}}
$$\tableau{1 & 2 \\ 1 & 3 \\ 2 & 4 & 3 \\ {\bf 2} & {\bf 4} & {\bf 3} & {\bf 1} & {\bf 5}}$$
\caption{A permuted basement SSAF with basement $(2,4,3,1,5)$ and shape $(3,3,1)$}
\end{figure}

\begin{proposition}{\label{keypermuted}}
The sum of the weights of all permuted basement SSAFs with partition shape $\lambda$ and basement $\omega$ is the key polynomial $\kappa_{\gamma}$, where $\gamma=\omega(\lambda)$.  (Symbolically, $\tilde{K}_{\omega,\lambda} = \kappa_{\omega(\lambda)}$.)
\end{proposition} 

\begin{proof}
Recall that the key polynomial $\kappa_{\gamma}$ consists of the weights of all SSAFs whose shape is less than or equal to $\gamma$ under the Bruhat ordering.  We construct a weight-preserving bijection between these SSAFs and the set $pb(\lambda, \omega)$, where $\omega(\lambda)=\gamma$.  In particular, this bijection preserves the row entries of each diagram.

Begin with an SSAF $F$ of shape less than or equal to $\gamma$ and an empty permuted basement SSAF, $\tilde{G}$, with basement $\omega$.  Find the largest entry, $\alpha$, in the first row of $F$ and place it on top of the leftmost entry (in the basement $\tilde{G}$) which is greater than or equal to $\alpha$.  Next find the second largest entry, $\beta$, in the first row of $F$ and place it on top of the leftmost available entry in the basement of $\tilde{G}$ that is greater than or equal to $\beta$.  (We say an entry is {\it available} if it lies beneath an empty cell.)  Continue in this manner until all of the entries in the first row of $F$ have been placed into $\tilde{G}$.  Repeat this process for each row (placing entries from row $i$ of $F$ into row $i$ of $\tilde{G}$) until all the row entries of $F$ have been inserted into $\tilde{G}$.  
 
 If $(a_1,a_2,\hdots, a_k)$ and $(b_1,b_2, \hdots, b_k)$ are two ordered sets of integers, we write $[a_1,a_2, \hdots, a_k] \le [b_1,b_2, \hdots, b_k]$ if $a_i \le b_i$ for all $i$ when the $a_i$ and $b_i$ are written in decreasing order.   If $F_i$ is the collection of entries in row $i$ of $F$ and $F_{i-1}'$ is the set containing the largest $|F_i|$ entries in row $i-1$ of $F$, then we must have $[F_i] \le [F_{i-1}']$.  Therefore the $j^{th}$ largest entry in row $i$ of $F$ must be less than or equal to at least $j$ entries in row $i-1$ of $F$, and hence after the placement of the largest $j-1$ entries of $F$ there is still a position in the $i^{th}$ row of $\tilde{G}$ in which the $j^{th}$ largest entry can be placed.  Therefore the map is well-defined. 
 
We must prove that the resulting diagram, $\tilde{G}$, is indeed a permuted basement SSAF.  To see this, we first show that $\tilde{G}$ satisfies the SSAF conditions.  The columns of $\tilde{G}$ are weakly increasing when read top to bottom by construction, so we must check that every triple is an inversion triple.  

Let the cells $\{a,b,c\}$ form a type $A$ triple in $\tilde{G}$ as shown.  $$\tableau{a & & b \\ c}$$  We know by construction that $a \le c$.  Since $b$ was not placed on top of $c$, either $b <a$ or $b > c$.  In either case, the three cells form a type $A$ inversion triple.

Next assume the cells $\{a,b,c\}$ form a type $B$ triple in $\tilde{G}$ as shown.  $$\tableau{ & & b \\ a & & c}$$  Then $b>a$ but $b \le c$.  So these cells form a type $B$ inversion triple.  Therefore every triple in $\tilde{G}$ is an inversion triple and hence $\tilde{G}$ satisfies the SSAF conditions.  

Now we prove that the shape of $\tilde{G}$ is a partition.  The shape of the SSAF $F$ must be less than or equal to $\gamma$.  This means that the columns of $F$ are obtained by permuting the columns of the partition $\lambda$ according to a permutation $\tau$ which is less than or equal to $\omega$ in the Bruhat order.  (From this point on we will use $\le$ to denote Bruhat inequality.)  The $i^{th}$ letter in the permutation $\tau$ determines which column of $F$ is the $i^{th}$ tallest column.  If the first row of $F$ contains $k_1$ non-empty cells, this set $S_1$ of cells are given by the first $k_1$ entries in $\tau$.  If $T_1$ is the collection of the first $k_1$ entries in $\omega$, then $[S_1] \le [T_1]$, since $\tau \le \omega$.  Therefore each of the entries in $S_1$ finds a position in the first $k_1$ columns of $\tilde{G}$.  

Assume that the second row of $F$ contains $k_2$ non-empty cells, collected into the set $F_2$.  These cells appear in the columns given by the first $k_2$ letters in $\tau$.  Therefore $[F_2] \le [S_2]$, where $S_2$ is the set consisting of the first $k_2$ letters of $\tau$.  Let $T_2$ be the set consisting of the first $k_2$ letters in $\omega$.  Then $[S_2] \le [T_2]$ since $\tau \le \omega$ in the Bruhat order.  The largest $k_2$ entries in row 2 of $F$ appear in the first $k_2$ columns of $\tilde{G}$, since $[F_2] \le [S_2] \le [T_2] \le [(T_1)]_{k_2}$, where $[(T_1)]_{k_2}$ is the set containing the $k_2$ largest entries in $T_1$.  Continuing this line of reasoning implies that for each row $i$ of $F$, the $k_i$ non-zero entries in row $i$ appear in the first $k_i$ columns of $\tilde{G}$.  Therefore the resulting shape is a partition shape.  So $\tilde{G}$ is indeed a permuted basement SSAF.  

Each SSAF in $\kappa_{\gamma}$ maps to a permuted basement SSAF $\tilde{G}$, so to prove that $\tilde{K}_{\gamma} = \kappa_{\gamma}$ we must show that these $\tilde{G}$ are the only permuted basement SSAFs with basement $\omega$.  To see this, we describe an inverse to the map above.  Consider an arbitrary permuted basement SSAF, $G$, with shape $\lambda$ and basement $\omega$.  Map the entries in each row of $G$ into an SSAF $F$ by the same process described above.  Again this mapping makes sense since the set $G_i$ of entries in the $i^{th}$ row of $G$ have the property that $[G_i] \le [G_{i-1}']$.  The result is an SSAF by the same argument as for $\tilde{G}$ above.

We must prove that the shape of the SSAF $F$ obtained from this mapping is less than or equal to $\omega$ in the Bruhat order on compositions.  We do this by proving that, for the permutation $\tau$ applied to $\lambda$ to obtain the shape of $F$, we have $[\tau(1), \tau(2), \hdots, \tau(k)] \le [\omega(1), \omega(2), \hdots, \omega(k)]$ for $1 \le k \le n$.  First notice that $\tau$ is given by listing the columns of $G$ in order from tallest to shortest, where columns of the same height are listed from left to right.  (See Figure 6.2 for an example.)

\begin{figure}
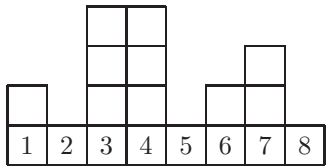
{\label{permword}}
$$\tableau{ & & {} & {} \\ & & {} & {} {} & & & {} \\ {} & & {} & {} & & {} & {} \\ 1 & 2 & 3 & 4 & 5 & 6 & 7 & 8}$$
\caption{The permutation associated to this composition is $\tau=(3,4,7,1,6,2,5,8)$.}
\end{figure}

Let $T_i$ be the set of entries in row $i$ of $G$ and let $k_i$ be the number of non-empty cells in row $i$ of $G$.  Then we have $[T_1] \le [\omega(1), \omega(2), \hdots, \omega(k_1)]$.  (For ease of notation, let $[\omega]_i=[\omega(1), \omega(2), \hdots, \omega(i)]$.)  The columns of $F$ containing the entries from $T_1$ represent the first $k_1$ letters of the permutation $\tau$ applied to $\lambda$ to obtain the shape of $F$.  Therefore $[\tau]_{k_1} \le  [\omega]_{k_1}$.  The remaining letters of $\tau$ are listed in increasing order following the first $k_1$ letters.  Let $m$ be the first position after the $k_1^{th}$ position where $\tau$ and $\omega$ differ.  Let $\alpha$ be the $m^{th}$ letter of $\tau$ and let $\beta$ be the $m^{th}$ letter of $\omega$.  Then $\alpha$ is the smallest positive integer that is not included in the first $m-1$ letters of $\tau$.  

We claim that $[\tau(1), \tau(2), \hdots, \tau(m-1), \alpha] \le [ \omega(1), \omega(2), \hdots \omega(m-1),\beta]$.  To see this, first note that $[\tau]_{m-1} \le [\omega]_{m-1}$, since these sets are obtained from $\{\tau(1),\tau(2), \hdots, \tau(k_1) \}$ and  $\{ \omega(1), \omega(2), \hdots \omega(k_1) \}$ by adding equivalent elements.  If $\beta > \alpha$, then $[\tau]_m \le [\omega]_m$ since we are adding a larger element to $[\omega]_{m-1}$ than to $[\tau]_{m-1}$.  Otherwise, let $\rho$ be the $(\alpha-1)^{th}$ largest element in $[\omega]_{m-1}$.  Then $\rho > \alpha-1$, since $\beta \le \alpha-1$ but $\beta \notin [\omega]_{m-1}$.  Adding $\beta$ to $[\omega]_{m-1}$ makes $\rho$ the $\alpha^{th}$ largest element in $[\omega]_m$ and therefore $\rho$ is greater than or equal to the $\alpha^{th}$ largest element, $\alpha$, in $[\tau]_m$.  The $\alpha-1$ smallest elements of $\tau$ are the set $\{ 1,2, \hdots, \alpha-1 \}$ and hence are less than any others set of positive integers consisting of $\alpha-1$ elements.  Therefore $[\tau]_m \le [\omega]_m$.  Repeat the above argument for each $m \ge k_1$ such that $\tau(m) \not= \omega(m)$.  This shows that $[\tau]_m \le [\omega]_m$ for all $m \ge k_1$.

Next let $r_2$ be the lowest row containing less than $k_1$ non-empty cells.  Then $[T_{r_2}] \le [T_{r_2-1}]$ and all of the entries in row $r_2$ of $F$ will appear in columns weakly to the left of the columns $\omega(1), \omega(2), \hdots, \omega(k_2)$, where $k_2$ is the number of entries in row $r_2$ of $G$.  Continuing this line of reasoning for each row of $G$, we see that the $i^{th}$ tallest column in $F$ appears weakly to the left of the column $\omega(i)$.  Since $\tau_i$ denotes the $i^{th}$ tallest column in $F$ (where columns of the same height are ordered from left to right), we have $[\tau(1), \tau(2) \hdots, \tau(j)] \le [\omega(1), \omega(2), \hdots, \omega(j)]$ for all $j \le k_1$.  Therefore the permutation applied to $\lambda$ to obtain the shape of $F$ is less than or equal to $\omega$ in the Bruhat order.
\end{proof}

Figure 6.3 depicts the key polynomial associated to the partition $(2,1)$ and permutation $(3,1,2)$.  Notice that when the basement of a permuted basement SSAF contains the permutation $\delta=(n,n-1, \hdots, 2, 1)$, the row entries appear in decreasing order.  Since every SSYT of shape $\lambda$ appears in $\kappa_{\delta(\lambda)}$, the collection of row-strict plane partitions with shape $\lambda$ is the set $pb(\lambda, \delta)$ with basements removed.  

One might be inclined to generalize the notion of a permuted basement SSAF to composition shapes.  The definition makes sense, but the collection of all such objects whose shape rearranges a fixed partition $\lambda$ is much larger than the collection of SSYT of shape $\lambda$.  An exploration of these notions and their connection to the Robinson-Schensted-Knuth algorithm will appear in \cite{HMR1}.

\begin{figure}
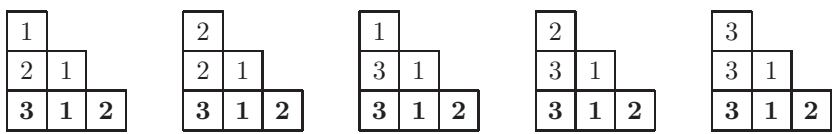

{\label{kappa}}
$$\tableau{1 \\ 2 & 1 \\ {\bf 3} & {\bf 1} & {\bf 2}}  \hspace*{.3in} \tableau{2 \\ 2 & 1 \\ {\bf 3} & {\bf 1} & {\bf 2}} \hspace*{.3in}  \tableau{1 \\ 3 & 1 \\ {\bf 3} & {\bf 1} & {\bf 2}}  \hspace*{.3in} \tableau{2 \\ 3 & 1 \\ {\bf 3} & {\bf 1} & {\bf 2}} \hspace*{.3in}  \tableau{3 \\ 3 & 1 \\ {\bf 3} & {\bf 1} & {\bf 2}}$$
\caption{$\kappa_{(3,1,2),(2,1)}=x_1^2x_2 + x_1x_2^2 + x_1^2x_3 + x_1 x_2 x_3 + x_1 x_3^2$}
\end{figure}

\bibliographystyle{amsalpha}

\end{document}